\documentclass[12pt]{article}
\usepackage{amsmath,amsfonts,amsthm,amssymb}
\newcommand{\F}{\mathbb F_q}
\newcommand{\K}{\overline{K}_c}
\newcommand{\AG}{\mathfrak A}
\newcommand{\MG}{\mathfrak M}
\newcommand{\FF}{\mathcal F}
\DeclareMathOperator{\gr}{gr}

\begin{document}
\newtheorem{lem}{Lemma}
\newtheorem{teo}{Theorem}
\newtheorem{prop}{Proposition}
\pagestyle{plain}
\title{Quasi-Holonomic Modules in Positive Characteristic}
\author{Anatoly N. Kochubei\footnote{Partially supported by
CRDF under Grant UM1-2567-OD-03, and by the Ukrainian Foundation 
for Fundamental Research, Grant 10.01/004.}\\ 
\footnotesize Institute of Mathematics,\\ 
\footnotesize National Academy of Sciences of Ukraine,\\ 
\footnotesize Tereshchenkivska 3, Kiev, 01601 Ukraine
\\ \footnotesize E-mail: \ kochubei@i.com.ua}
\date{}
\maketitle
\newpage
\begin{abstract}
We study modules over the Carlitz ring, a counterpart
of the Weyl algebra in analysis over local fields of positive 
characteristic. It is shown that some basic objects of function
field arithmetic, like the Carlitz module, Thakur's
hypergeometric polynomials, and analogs of binomial coefficients 
arising in the function field version of umbral calculus, 
generate quasi-holonomic modules. This class of modules is, in 
many respects, similar to the class of holonomic modules in the 
characteristic zero theory.
\end{abstract}
\vspace{2cm}
\newpage

\section{INTRODUCTION}

The theory of holonomic modules over the Weyl algebra and more 
general algebras of differential or $q$-difference operators is 
becoming increasingly important, both as a crucial part of the 
general theory of D-modules and in view of various applications 
(see, for example, \cite{BK,Cart,Gue,S}). Well-known pathological 
properties of differential operators over fields of positive 
characteristic make the available, for this case, analogs of the 
theory of D-modules much more complicated \cite{Bog,Lyub}. More 
importantly, the resulting structures are not connected with the 
existing analysis in positive characteristic based on a completely 
different algebraic foundation.

Any non-discrete locally compact field of a positive 
characteristic $p$ is isomorphic to the field $K$ of formal Laurent 
series with coefficients from the Galois field $\F$, $q=p^\nu$, $\nu \in \mathbb 
Z_+$. The field $K$ is endowed with a non-Archimedean absolute 
value as follows. If $z\in K$,
$$
z=\sum\limits_{i=m}^\infty \zeta_ix^i,\quad m\in \mathbb 
Z,\ \zeta_i\in \F ,\ \zeta_m\ne 0,
$$
then $|z|=q^{-m}$. This valuation can be extended onto the field 
$\K$, the completion of an algebraic closure of $K$.

Analysis over $K$ and $\K$, which was initiated in the great paper 
by Carlitz \cite{Carl} and developed subsequently by Wagner, Goss, 
Thakur, the author, and many others (see the bibliography in 
\cite{G2,Th3}) is very different from the classical calculus. An 
important feature is the availability of many non-trivial additive 
(actually, $\F$-linear) polynomials and power series of the form 
$u(t)=\sum\limits a_kt^{q^k}$.

Taking into account the fact that the usual factorial $i!$, seen as 
an element of $K$, vanishes for $i\ge p$, Carlitz introduced the new factorial
\begin{equation} 
D_i=[i][i-1]^q\ldots [1]^{q^{i-1}},\quad [i]=x^{q^i}-x\ (i\ge 1),\ 
D_0=1,
\end{equation}
the $\F$-linear logarithm and exponential (which obtained a wide 
generalization later, in the theory of Drinfeld modules), as well 
as an important polynomial system, the Carlitz polynomials. 
Subsequently many other $\F$-linear special functions, such as 
Thakur's hypergeometric function \cite{Th1,Th2,Th3} and further special
polynomial systems, were introduced and investigated. The 
difference operator
\begin{equation}
\Delta u(t)=u(xt)-xu(t)
\end{equation} 
introduced in \cite{Carl} became the main ingredient of the 
$\F$-linear calculus and analytic theory of differential equations 
over $K$ developed in \cite{K2,K3,K4}. The role of a derivative is 
played by the $\F$-linear operator $d=\sqrt[q]{}\circ \Delta$ ({\it the 
Carlitz derivative}). The latter appears also in the 
$\F$-linear umbral calculus \cite{K5} where an important role 
belongs to the following new analog of binomial coefficients
\begin{equation} 
\binom{k}{m}_K=\frac{D_k}{D_mD_{k-m}^{q^m}},\quad 0\le m\le k.
\end{equation}
The meaning of a polynomial coefficient in a differential equation 
of the above type is not a usual multiplication by a polynomial, 
but the action of a polynomial in the Frobenius operator $\tau$, 
$\tau u=u^q$. With this notation, $d=\tau^{-1}\Delta$. The 
operator $d$ is defined on any $\F$-linear $\K$-valued continuous 
function; in particular, it decreases by one the ``$\F$-linear 
degree'' of any $\F$-linear polynomial (see the relation (8) 
below).

The above developments show that in the positive characteristic 
case a natural counterpart of the Weyl algebra is, for the case of 
a single variable, the ring $\AG_1$ generated by $\tau,d$, and 
scalars from $\K$, with the relations \cite{K1}
\begin{equation}
d\tau -\tau d=[1]^{1/q},\quad \tau \lambda =\lambda^q\tau ,\quad ,
d\lambda =\lambda^{1/q}d\ (\lambda \in \K).
\end{equation} 
Some algebraic properties of $\AG_1$ were studied in \cite{K3} -- 
it is left and right Noetherian, with no zero divisors.

The aim of this paper is to initiate the dimension theory for 
modules over $\AG_1$ and more general ``several variable'' rings. 
The definition of the latter is not straightforward. If, for 
example, we consider the natural action of the Carlitz derivatives 
$d_s$ and $d_t$ on an $\F$-linear monomial 
$f(s,t)=s^{q^m}t^{q^n}$, we notice immediately that $d_s^mf$ is 
not a polynomial, nor even a holomorphic function in $t$, if 
$m>n$ (since the action of $d$ is not linear and involves taking 
the $q$-th root). Moreover, it follows from the relation 
$d(s^{q^m})=[m]^{1/q}s^{q^{m-1}}$ and the last commutation 
relation in (4) that $d_s$ and $d_t$ do not commute even on 
monomials $f$ with $m<n$.

A reasonable generalization is inspired by Zeilberger's idea (see 
\cite{Cart}) to study holonomic properties of sequences of 
functions making a transform with respect to the discrete 
variables, which reduces the continuous-discrete case to the 
purely continuous one (simultaneously in all the variables). In 
our situation, if $\{ P_k(s)\}$ is a sequence of $\F$-linear 
polynomials with $\deg P_k\le q^k$, we set
\begin{equation*}
f(s,t)=\sum \limits_{k=0}^\infty P_k(s)t^{q^k},\tag{$*$}
\end{equation*}
and $d_s$ is well-defined. In the variable $t$, we consider not 
$d_t$ but the linear operator $\Delta_t$. The latter does not 
commute with $d_s$ either, but satisfies the commutation relations
$$
d_s\Delta_t-\Delta_td_s=[1]^{1/q}d_s,\quad \Delta_t\tau -\tau 
\Delta_t=[1]\tau ,
$$
so that the resulting ring $\AG_2$ resembles a universal 
enveloping algebra of a solvable Lie algebra. Similarly we define 
$\AG_{n+1}$ for $n>1$.

Introducing in $\AG_{n+1}$ an analog of the Bernstein filtration 
and considering filtered modules over $\AG_{n+1}$, we find that 
basic principles of the theory of algebraic D-modules \cite{Cout} 
carry over to this case without serious complications. However, 
the nonlinearity of $\tau$ and $d$ brings new phenomena. In 
particular, already the ring $\AG_1$ possesses non-trivial 
finite-dimensional representations. Therefore an analog of the 
Bernstein inequality does not hold here without some additional 
assumptions.

In spite of this fact, the notion of a holonomic module (that is a 
module with the minimal possible GK dimension) seems to have a 
reasonable sense for the case of $\AG_{n+1}$-modules. The examples 
considered in this paper (both for $\AG_1$-modules and 
$\AG_{n+1}$-modules with $n\ge 1$) show that the cases of an 
anomalously small GK dimension may be seen as degenerate ones. In 
terms of applications to analysis, it appears that a remarkable 
phenomenon discovered by Zeilberger (see \cite{Cart}) -- that 
virtually all important special functions and sequences of 
classical analysis generate holonomic modules -- is maintained in 
the positive characteristic case, if a holonomic module is defined 
as a one with a minimal ``generic'' GK dimension, with degenerate 
cases excluded. In the author's opinion, such applications provide 
a sufficient justification for the definition of a quasi-holonomic 
module given in this paper (Sect. 3.2).

Accordingly, the case we study in a greater detail is that of quasi-holonomic 
submodules of the $\AG_{n+1}$-module of $\F$-linear functions 
$u(s,t_1,\ldots ,t_n)$, polynomial in $s$ and holomorphic near the 
origin in $t_1,\ldots ,t_n$. Following \cite{Cart} we call a 
function $f$ quasi-holonomic if such is the module $\AG_{n+1}f$. We 
prove general conditions for a function $f$ to be quasi-holonomic and 
verify them for basic objects of this branch of analysis -- the 
Carlitz polynomials, Thakur's hypergeometric polynomials, and the 
$K$-binomial coefficients (3), making the above transformation (*) from 
discrete variables to continuous ones.

Considering the $K$-binomial coefficients we use this occasion to 
prove also the fact that they belong to the ring of integers not 
only for the field $K$, but for any place of the global function 
field $\F (x)$. Together with the results of \cite{K5}, this 
property supports the case for considering the expressions (3) as 
``proper'' analogs of the classical binomial coefficients. For 
other analogs of the latter see \cite{Th3}.

\section{The Carlitz Ring}

{\bf 2.1.} Denote by $\FF_{n+1}$ the set of all germs of functions of the 
form
\begin{equation}
f(s,t_1,\ldots ,t_n)=\sum\limits_{k_1=0}^\infty \ldots \sum\limits_{k_n=0}^\infty
\sum\limits_{m=0}^{\min (k_1,\ldots ,k_n)}a_{m,k_1,\ldots 
,k_n}s^{q^m}t_1^{q^{k_1}}\ldots t_n^{q^{k_n}}
\end{equation}
where $a_{m,k_1,\ldots ,k_n}\in \K$ are such that all the series 
are convergent on some neighbourhoods of the origin. We do not 
exclude the case $n=0$ where $\FF_1$ will mean the set of all 
$\F$-linear power series $\sum\limits_ma_ms^{q^m}$ convergent on 
a neighbourhood of the origin. $\widehat{\FF}_{n+1}$ will denote 
the set of all polynomials from $\FF_{n+1}$, that is the series 
(5) in which only a finite number of coefficients is different 
from zero.

The ring $\AG_{n+1}$ is generated by the operators $\tau 
,d_s,\Delta_{t_1},\ldots \Delta_{t_n}$ on $\FF_{n+1}$ defined in 
the Introduction, and the operators of multiplication by scalars 
from $\K$. To simplify the notation, we will write $\Delta_j$ 
instead of $\Delta_{t_j}$ and identify a scalar $\lambda \in \K$ 
with the operator of multiplication by $\lambda$. The operators 
$\Delta_j$ are $\K$-linear, so that
\begin{equation}
\Delta_j\lambda =\lambda \Delta_j,\quad \lambda \in \K ,
\end{equation}
while the operators $\tau ,d_s$ satisfy the commutation relations 
(4). In the action of each operator $d_s,\Delta_j$ (acting in a 
single variable), other variables are treated as scalars. The 
operator $\tau$ acts simultaneously on all the variables and 
coefficients, so that
$$
\tau f=\sum a^q_{m,k_1,\ldots 
,k_n}s^{q^{m+1}}t_1^{q^{k_1+1}}\ldots t_n^{q^{k_n+1}}.
$$

It follows from (2) that 
\begin{equation}
\Delta_jt_j^{q^k}=\begin{cases}
[k]t_j^{q^k}, & \text{if $k\ge 1$};\\
0, & \text{if $k=0$};
\end{cases}
\end{equation}
the second equality can be included in the first one, if we set $[0]=0$.
Similarly
\begin{equation}
d_ss^{q^m}=[m]^{1/q}s^{q^{m-1}},\quad m\ge 0.
\end{equation} 
Since $|[m]|=q^{-1}$ for any $m\ge 1$, the action of operators 
from $\AG_{n+1}$ does not spoil convergence of the series (5).

The identity $[k+1]-[k]^q=[1]$, together with (7) and (8), implies 
the commutation relations
\begin{equation}
\Delta_j\tau -\tau \Delta_j=[1]\tau ,\quad 
d_s\Delta_j-\Delta_jd_s=[1]^{1/q}d_s,\quad j=1,\ldots ,n,
\end{equation} 
verified by applying both sides of each equality to an arbitrary 
monomial.

Using the commutation relations (4), (6), and (9), we can write 
any element $a\in \AG_{n+1}$ as a finite sum
\begin{equation}
a=\sum c_{l,\mu,i_1,\ldots ,i_n}\tau^ld_s^\mu \Delta_1^{i_1}\ldots 
\Delta_n^{i_n}.
\end{equation} 

\medskip
\begin{prop}
The representation (10) of an  element $a\in \AG_{n+1}$ is unique.
\end{prop}

\medskip
{\it Proof}. Suppose that 
\begin{equation}
\sum\limits_{l,\mu,i_1,\ldots ,i_n}c_{l,\mu,i_1,\ldots ,i_n}\tau^ld_s^\mu 
\Delta_1^{i_1}\ldots \Delta_n^{i_n}=0.
\end{equation} 
Applying the left-hand side of (11) to the function 
$st_1^{q^{k_1}}\ldots t_n^{q^{k_n}}$ with $k_1,\ldots ,k_n>0$ we 
find that
$$
\sum\limits_l\left( \sum\limits_{i_1,\ldots ,i_n}c_{l,0,i_1,\ldots ,i_n}
[k_1]^{i_1q^l}\ldots [k_n]^{i_nq^l}\right) 
s^{q^l}t_1^{q^{k_1+l}}\ldots t_n^{q^{k_n+l}}=0
$$
whence
$$
\sum\limits_{i_1,\ldots ,i_n}c_{l,0,i_1,\ldots ,i_n}
[k_1]^{i_1q^l}\ldots [k_n]^{i_nq^l}=0
$$
for each $l$. Writing this in the form
\begin{equation}
\sum\limits_{i_n}\rho (i_n)y^{i_n}=0
\end{equation} 
where
$$
\rho (i_n)=\sum\limits_{i_1,\ldots ,i_{n-1}}c_{l,0,i_1,\ldots ,i_n}
[k_1]^{i_1q^l}\ldots [k_{n-1}]^{i_{n-1}q^l},\quad y=[k_n]^{q^l},
$$
and taking into account that (12) holds for arbitrary $k_n\ge 1$, 
that is for an infinite set of values of $y$, we find that $\rho 
(i_n)=0$. Repeating this reasoning we get the equality $c_{l,0,i_1,\ldots 
,i_n}=0$ for all $l,0,i_1,\ldots ,i_n$.

Suppose that $c_{l,\mu ,i_1,\ldots ,i_n}=0$ for $\mu \le \mu_0$ 
and arbitrary $l,i_1,\ldots ,i_n$. Then we apply the left-hand 
side of (11) to the function $s^{q^{\mu_0+1}}t_1^{q^{k_1}}\ldots 
t_n^{q^{k_n}}$ and proceed as before coming to the equality 
$c_{l,\mu_0+1,i_1,\ldots ,i_n}=0$ for all $l,i_1,\ldots ,i_n$. 
$\qquad \blacksquare$

\medskip
It is easy to prove by induction with respect to $n$ (using the 
commutation relations (9) and the result from \cite{K3} regarding 
the case $n=0$) that $\AG_{n+1}$ has no zero-divisors.

\medskip
{\bf 2.2.} Let us introduce a filtration in $\AG_{n+1}$ denoting 
by $\Gamma_\nu$, $\nu \in \mathbb Z_+$, the $\K$-vector space of 
operators (10) with $\max \{ l+\mu +i_1+\cdots +i_n\} \le \nu$ 
where the maximum is taken over all the terms contained in the 
representation (10). It is clear that $\AG_{n+1}$ is a filtered 
ring (for the definitions see \cite{MR}). Setting $T_0=\K$, $T_\nu 
=\Gamma_\nu /\Gamma_{\nu -1}$, $\nu \ge 1$, we introduce the 
associated graded ring
$$
\gr (\AG_{n+1})=\bigoplus\limits_{\nu=0}^\infty T_\nu .
$$
It is generated by scalars $\lambda \in T_0$ and the images 
$\bar \tau,\bar d_s,\bar \Delta_1,\ldots 
,\bar \Delta_n\in T_1$ of the elements $\tau 
,d_s,\Delta_1,\ldots ,\Delta_n\in \Gamma_1$ respectively, which 
satisfy, by virtue of (4), (6), and (9), the relations
\begin{gather*}
\bar d_s\bar\tau-\bar\tau \bar d_s=0,
\bar\tau \lambda 
=\lambda^q\bar\tau,\bar d_s\lambda 
=\lambda^{1/q}\bar d_s,\\
\bar d_s\bar \Delta_j-\bar \Delta_j\bar d_s=0,
\bar \Delta_j\bar \tau -\bar\tau \bar \Delta_j=0,
\bar \Delta_j\lambda =\lambda \bar \Delta_j\quad 
(j=1,\ldots ,n).
\end{gather*}

It is clear that $\AG_{n+1}$ is a (left and right) almost 
normalizing extension of the field $\K$ (see Chapter 1, \S 6 in 
\cite{MR}), so that the rings $\AG_{n+1}$ and $\gr (\AG_{n+1})$ 
are left and right Noetherian.

Let us compute the dimension of the $\K$-vector space 
$\Gamma_\nu$. Note that
$$
\dim \Gamma_\nu =\dim \bigoplus\limits_{j=1}^\nu T_j,
$$
so that $\dim \Gamma_\nu$ coincides with the dimension of the 
appropriate space appearing in the natural filtration in $\gr (\AG_{n+1})$.

\medskip
\begin{lem}
For any $\nu \in \mathbb N$
$$
\dim \Gamma_\nu = \binom{\nu +n+2}{n+2}.
$$
\end{lem}

\medskip
{\it Proof}. The number $\dim \Gamma_\nu$ coincides with the 
number of non-negative integral solutions $(l,\mu ,i_1,\ldots 
,i_n)$ of the inequality $l+\mu +i_1+\cdots +i_n\le \nu$, so that
$$
\dim \Gamma_\nu =\sum\limits_{j=0}^\nu N(j,n+2)
$$
where $N(j,k)$ is the number of different representations of $j$ 
as sums of $k$ non-negative integers. It is known (Proposition 6.1 in
\cite{Lan}) that $N(j,k)=\dbinom{j+k-1}{k-1}$. Then (see Sect. 1.3 
from \cite{Ri})
$$
\dim \Gamma_\nu =\sum\limits_{j=0}^\nu \binom{j+n+1}{n+1}=
\sum\limits_{i=0}^\nu \binom{\nu+n+1-i}{n+1}=\binom{\nu 
+n+2}{n+2},
$$
as desired. $\qquad \blacksquare$

\section{Filtered Modules}

{\bf 3.1.} Let $M$ be a left module over the Carlitz ring 
$\AG_{n+1}$. Suppose we have a filtration $\{ \MG_j\}$ of $M$, 
that is
\begin{equation}
\MG_0\subset \MG_1\subset \ldots \subset M,\quad 
M=\bigcup\limits_{j\ge 0}\MG_j,
\end{equation}
and $\Gamma_\nu \MG_j\subset \MG_{\nu +j}$ for any $\nu ,j\in 
\mathbb Z_+$. We assume that each $\MG_j$ is a finite-dimensional 
vector space over $\K$. Below we write $\MG_j=\{ 0\}$ and 
$\Gamma_\nu =\{ 0\}$ if $j<0$ and $\nu <0$.

In a standard way \cite{Cout} we define the graded module
$$
\gr (M)=\bigoplus\limits_{j\ge 0}\left( \MG_j/\MG_{j-1}\right)
$$
over $\gr (\AG_{n+1})$, associated with the filtration (13). As 
usual, the filtration (13) is called {\it good}, if $\gr (M)$ is 
finitely generated.

Main properties of filtered modules over the Weyl algebra (see 
\cite{Bj,Cout}) carry over to our situation without any 
substantial changes, both in their formulations and proofs. In 
fact, the only technical difference is that the operators $\tau$ 
and $d_s$ are semilinear, not linear. However, as it is explained 
in Appendix I to Chapter 2 of \cite{Bour}, basic notions of linear 
algebra remain valid for semilinear mappings -- a semilinear 
mapping of a vector space into itself can be interpreted as a 
linear mapping between two different vector spaces, and, for 
instance, dimensions of the kernel and cokernel are not changed in 
this interpretation. Note that everywhere in this paper we 
consider vector spaces over the algebraically closed field $\K$, 
on which $\tau$ induces an automorphism. Below, as before, $\dim$ 
means the dimension over $\K$.

In particular, for a good filtration there exist a polynomial 
$\chi \in \mathbb Q[t]$ and a number $N\in \mathbb N$, such that
$$
\dim \MG_s=\sum\limits_{i=0}^s\dim (\MG_i/\MG_{i-1})=\chi 
(s)\text{ for }s\ge N.
$$
The number $d(M)=\deg \chi$, called the (Gelfand-Kirillov) {\it 
dimension} of $M$, and the leading coefficient of $\chi$ 
multiplied by $d(M)!$, called the {\it multiplicity} $m(M)$ of 
$M$, do not depend on the choice of a good filtration on $M$. A 
filtration $\{ \MG_i\}$ is good if and only if there exists such 
$k_0\in \mathbb N$ that
$$
\MG_{i+k}=\Gamma_i\MG_k\text{ \ for all }k\ge k_0.
$$

If $N$ and $M/N$ are a submodule and the corresponding quotient 
module, with the induced filtrations, then $d(M)=\max \{ 
d(N),d(M/N)\}$, and if $d(N)=d(M/N)$, then $m(M)=m(N)+m(M/N)$. For 
a direct sum $M=M_1\oplus \cdots \oplus M_k$ we have $d(M)=\max \{ 
d(M_1),\ldots ,d(M_k)\}$.

In particular, if we consider $\AG_{n+1}$ as a left module over 
itself, then by Lemma 1
\begin{equation}
d(\AG_{n+1})=n+2,\quad m(\AG_{n+1})=1.
\end{equation}
It follows from (14) and the above general facts that for any 
finitely generated left $\AG_{n+1}$-module
\begin{equation}
d(M)\le n+2.
\end{equation}

By (14), the bound in (15) in general cannot be improved. However, 
if $I$ is a non-zero left ideal in $\AG_{n+1}$, then
\begin{equation}
d(\AG_{n+1}/I)\le n+1.
\end{equation}
The proof of (16) is identical to the proof of Corollary 9.3.5 
from \cite{Cout}.

\medskip
{\bf 3.2.} Let us consider the set $\widehat{\FF}_{n+1}$ of 
polynomials (5) as a $\AG_{n+1}$-module. A filtration
$$
\FF^{(0)}_{n+1}\subset \FF^{(1)}_{n+1}\subset \ldots \subset 
\widehat{\FF}_{n+1}
$$
can be introduced by setting $\FF^{(j)}_{n+1}$ to be the 
collection of all the polynomials (5), in which the maximal 
indices $k_1,\ldots ,k_n$ corresponding to non-zero coefficients 
$a_{m,k_1,\ldots ,k_n}$ do not exceed $j$. This filtration is 
obviously good.

\medskip
\begin{prop}
For the module $\widehat{\FF}_{n+1}$,
\begin{equation}
d\left( \widehat{\FF}_{n+1}\right) =n+1,\quad m\left( \widehat{\FF}_{n+1}\right) 
=n!
\end{equation}
\end{prop}

\medskip
{\it Proof}. Let us compute $\dim \FF^{(j)}_{n+1}$. For a fixed 
$\mu$, the quantity of $n$-tuples $(k_1,\ldots ,k_n)$ of 
non-negative integers, for which $\min (k_1,\ldots ,k_n)=\mu$, is 
added up from those $n$-tuples where $i$ numbers are equal to 
$\mu$ while $n-i$ numbers are strictly larger and can take $j-\mu$ 
values. Therefore the above quantity equals 
$\sum\limits_{i=1}^n\dbinom{n}{i}(j-\mu )^{n-i}$. Next, $\mu +1$ 
possible values of $m$ in (5) correspond to each $n$-tuple. Thus,
$$
\dim \FF^{(j)}_{n+1}=\sum\limits_{\mu =0}^j(\mu +1)
\sum\limits_{i=1}^n\dbinom{n}{i}(j-\mu )^{n-i}=\sum\limits_{\mu =0}^j(\mu 
+1)\left\{ (j-\mu +1)^n-(j-\mu )^n\right\} .
$$

Denote $r_\mu =(j-\mu +1)^n-(j-\mu )^n$, $R_i=r_0+r_1+\cdots 
+r_i=(j+1)^n-(j-i)^n$. Performing the Abel transformation we get
\begin{multline*}
\dim \FF^{(j)}_{n+1}=(j+1)R_j-\sum\limits_{i=0}^{j-1}R_i
=(j+1)^{n+1}-j(j+1)^n+\sum\limits_{i=0}^{j-1}(j-i)^n\\
=(j+1)^n+\sum\limits_{k=1}^jk^n=(j+1)^n+S_n(j+1)
\end{multline*}
where $S_n(N)=1^n+2^n+\cdots +(N-1)^n$.

It is known (\cite{IR}, Chapter 15) that
$$
S_n(N)=\frac{1}{n+1}\sum\limits_{k=0}^n\binom{n+1}{k}B_kN^{n+1-k}
$$
where $B_k$ are the Bernoulli numbers. Therefore we find that
$$
\dim \FF^{(j)}_{n+1}=\frac{(j+1)^{n+1}}{n+1}+P_n(j)
$$
where $P_n$ is a polynomial of the degree $n$. This implies (17). 
$\qquad \blacksquare$

\medskip
It is natural to call an $\AG_{n+1}$-module $M$ {\it quasi-holonomic} if 
$d(M)=n+1$. Thus, $\widehat{\FF}_{n+1}$ is an example of a 
quasi-holonomic module.

\medskip
{\bf 3.3.} Let us look at possible values of $d(M)$ for 
$\AG_1$-modules. The next result demonstrates a sharp difference 
from the case of modules over the Weyl algebras.

\medskip
\begin{teo}
\begin{description}
\item{{\rm (i)}} For any $k=1,2,\ldots$, there exists such a 
nontrivial $\AG_1$-module $M$ that $\dim M=k$ ($\dim$ means the 
dimension over $\K$), that is $d(M)=0$.
\item{{\rm (ii)}} Let $M$ be a finitely generated $\AG_1$-module 
with a good filtration. Suppose that there exists a ``vacuum vector'' $v\in 
M$, such that $d_sv=0$ and $\tau^m(v)\ne 0$ for all 
$m=0,1,2,\ldots$. Then $d(M)\ge 1$.
\end{description}
\end{teo}

\medskip
{\it Proof}. (i) Let $M=(\K )^k$. Denote by $\mathbf e_1,\ldots 
,\mathbf e_k$ the standard basis in $M$, that is $\mathbf 
e_j=(0,\ldots ,0,1,0,\ldots ,0)$, with 1 at the $j$-th place. Let 
$(\lambda_{ij})$ be a $k\times k$ matrix over $\K$, such that 
$\lambda_{ij}\in \F$ if $i\ne j$, while the diagonal elements 
satisfy the equation $\lambda^q-\lambda +[1]^{1/q}=0$. We define 
the action of $\tau$ and $d_s$ on $M$ as follows:
$$
\tau (c\mathbf e_j)=c^q\mathbf e_j;\ d_s(\mathbf 
e_j)=\sum\limits_{i=1}^n\lambda_{ij}\mathbf e_i;\ d_s(c\mathbf 
e_j)=c^{1/q}\mathbf e_j,\quad c\in \K ,j=1,\ldots ,k,
$$
with subsequent additive continuation onto $M$. 

If $x=\sum\limits_{j=1}^kc_j\mathbf e_j$, $c_j\in \K$, then we 
have
$$
\tau d_s(x)=\sum\limits_{j=1}^kc_j\sum\limits_{i=1}^n\lambda_{ij}^q\mathbf 
e_i,\quad d_s\tau (x)=\sum\limits_{j=1}^kc_j\sum\limits_{i=1}^n\lambda_{ij}\mathbf 
e_i,
$$
so that
$$
d_s\tau (x)-\tau d_s(x)=[1]^{1/q}x,
$$
and we have indeed an $\AG_1$-module.

(ii) It follows from the relation 
$[d_s,\tau^m]=[m]^{1/q}\tau^{m-1}$ (see \cite{K3}) that
$$
d_s\tau^mv=[m]^{1/q}\tau^{m-1}v,\quad m=1,2,\ldots ,
$$
that is $\tau^{m-1}v$ is an eigenvector of a linear operator 
$d_s\tau$ on $M$ (considered as a $\K$-vector space) corresponding 
to the eigenvalue $[m]^{1/q}$. Therefore the vectors $\tau^{m-1}v$ 
are linearly independent. It follows from the existence of the 
Hilbert polynomial $\chi$ implementing the dimension $d(M)$ that 
$d(M)\ge 1$. $\qquad \blacksquare$

\medskip
\section{Holonomic Functions}

{\bf 4.1.} Let $0\ne f\in \FF_{n+1}$,
$$
I_f=\left\{ \varphi \in \AG_{n+1}:\ \varphi (f)=0\right\} .
$$
$I_f$ is a left ideal in $\AG_{n+1}$. The left $\AG_{n+1}$-module 
$M_f=\AG_{n+1}/I_f$ is isomorphic to the submodule 
$\AG_{n+1}f\subset \FF_{n+1}$ -- an element $\varphi (f)\in 
\AG_{n+1}f$ corresponds to the class of $\varphi \in \AG_{n+1}$ in 
$M_f$. A natural good filtration in $M_f$ is induced from that in 
$\AG_{n+1}$ -- the subspace $\MG_j$ is generated by elements 
$\tau^ld_s^\mu \Delta_1^{i_1}\ldots \Delta_n^{i_n}f$ with $l+\mu 
+i_1+\cdots +i_n\le j$.

As we know (see (16)), if $I_f\ne \{ 0\}$, then $d(M_f)\le n+1$. 
We call a function $f$ {\it quasi-holonomic} if the module $M_f$ is 
quasi-holonomic, that is $d(M_f)=n+1$. The condition $I_f\ne \{ 0\}$ 
means that $f$ is a solution of a ``differential equation'' 
$\varphi (f)=0$, $\varphi \in \AG_{n+1}$. For $n=0$, we have the 
following easy result.

\medskip
\begin{teo}
If a non-zero function $f\in \FF_1$ satisfies an equation
$\varphi (f)=0$, $0\ne \varphi \in \AG_1$, then $f$ is quasi-holonomic.
\end{teo}

\medskip
{\it Proof}. It is sufficient to show that $\dim M_f=\infty$. In 
fact, the sequence $\left\{ \tau^lf\right\}_{l=0}^\infty$ is 
linearly independent because otherwise we would have such a 
finite collection of elements $c_0,c_1,\ldots ,c_N\in \K$, some of 
which are different from zero, that
\begin{equation}
c_0f(s)+c_1f^q(s)+\cdots +c_Nf^{q^N}(s)=0
\end{equation}
for all $s$ from a neighbourhood of the origin in $\K$. It follows 
from (18) that $f$ takes only a finite number of values. By the 
uniqueness theorem for non-Archimedean holomorphic functions, 
$f(s)\equiv \text{const}$ on some neighbourhood of the origin. Due 
to the $\F$-linearity, $f(s)\equiv 0$, and we have come to a 
contradiction. $\qquad \blacksquare$

\medskip
In particular, any $\F$-linear polynomial of $s$ is quasi-holonomic, 
since it is annihilated by $d_s^m$, with a sufficiently large $m$.

\medskip
{\bf 4.2.} If $n>0$, the situation is more complicated. We call 
the module $M_f$ (and the corresponding function $f$) {\it 
degenerate} if $d(M_f)<n+1$ (by the Bernstein inequality, there is 
no degeneracy phenomena for modules over the complex Weyl 
algebra). We give an example of degeneracy for the case $n=1$.

Let $f(s,t_1)=g(st_1)\in \FF_2$ where the function $g$ belongs to
$\FF_1$ and satisfies an equation $\varphi (g)=0$, $\varphi \in 
\AG_1$. Then $f$ is degenerate.

Indeed, by the general rule, $\MG_j$ is spanned by elements 
$\tau^ld_s^\mu \Delta_1^{i_i}f$ with $l+\mu +i_1\le j$. In the present 
situation,  
$$
\Delta_1f=g(xst_1)-xg(st_1)=\tau d_sg,
$$
so that an element $\tau^ld_s^\mu \Delta_1^{i_i}f$ is a linear 
combination of elements $\left( \tau^{l+\lambda}d_s^{\mu 
+\nu}g\right) (s,t)$ with $\lambda \le i_1$, $\nu \le i_1$. 
Therefore $\MG_j$ is contained in the linear hull of elements 
$\tau^kd_s^mg$, $k+m\le 2j$. By Theorem 2, the $\K$-dimension of 
the latter does not exceed a linear function of $2j$, so that 
$d(M_f)\le 1$. On the other hand, since, as in the proof of 
Theorem 2, the system of functions $\left\{ 
\tau^lf\right\}_{l=0}^\infty$ is linearly independent, we find that
$d(M_f)=1$.

In order to exclude the degenerate case, we introduce the notion 
of a non-sparse function.

A function $f\in \FF_{n+1}$ of the form (5) is called {\it 
non-sparse} if there exists such a sequence $m_l\to \infty$ that, 
for any $l$, there exist sequences $k_1^{(i)},k_2^{(i)},\ldots 
,k_n^{(i)}\ge m_l$ (depending on $l$), such that $k_\nu^{(i)}\to 
\infty$ as $i\to \infty$ ($\nu =1,\ldots ,n$), and $a_{m,k_1^{(i)},\ldots 
,k_n^{(i)}}\ne 0$.

\medskip
\begin{lem}
If a function $f$ is non-sparse, then the system of functions 
$(\tau d_s)^\lambda \Delta_1^{j_1}\ldots \Delta_n^{j_n}f$ 
($\lambda ,j_1,\ldots ,j_n=0,1,2,\ldots$) is linearly independent 
over $\K$.
\end{lem}

\medskip
{\it Proof}. Suppose that
\begin{equation}
\sum\limits_{\lambda =0}^\Lambda \sum\limits_{j_1=0}^{J_1}\ldots 
\sum\limits_{j_n=0}^{J_n}c_{\lambda ,j_1,\ldots ,j_n}
(\tau d_s)^\lambda \Delta_1^{j_1}\ldots \Delta_n^{j_n}f=0
\end{equation}
for some $c_{\lambda ,j_1,\ldots ,j_n}\in \K$, $\Lambda ,J_1,\ldots 
,J_n\in \mathbb N$. Substituting (5) into (19) and collecting 
coefficients of the power series we find that
\begin{equation}
\sum\limits_{\lambda =0}^\Lambda \sum\limits_{j_1=0}^{J_1}\ldots 
\sum\limits_{j_n=0}^{J_n}c_{\lambda ,j_1,\ldots ,j_n}
[m_l]^\lambda [k_1^{(i)}]^{j_1}\ldots[k_n^{(i)}]^{j_n}=0
\end{equation}
for all $l,i$.

We see from (20) that the polynomial
$$
\sum\limits_{j_n=0}^{J_n}\left\{ \sum\limits_{\lambda =0}^\Lambda 
\sum\limits_{j_1=0}^{J_1}\ldots 
\sum\limits_{j_{n-1}=0}^{J_{n-1}}c_{\lambda ,j_1,\ldots ,j_n}
[m_l]^\lambda [k_1^{(i)}]^{j_1}\ldots[k_{n-1}^{(i)}]^{j_{n-1}}
\right\} z^{j_n}
$$
has an infinite sequence of different roots, so that
$$
\sum\limits_{\lambda =0}^\Lambda 
\sum\limits_{j_1=0}^{J_1}\ldots 
\sum\limits_{j_{n-1}=0}^{J_{n-1}}c_{\lambda ,j_1,\ldots ,j_n}
[m_l]^\lambda [k_1^{(i)}]^{j_1}\ldots[k_{n-1}^{(i)}]^{j_{n-1}}=0
$$
for all $l,i$, and for each $j_n=0,1,\ldots ,J_n$. Repeating this 
reasoning we find that all the coefficients $c_{\lambda ,j_1,\ldots 
,j_n}$ are equal to zero. $\qquad\blacksquare$

\medskip
Now the above arguments regarding $d(M_f)$ yield the following 
result.

\medskip
\begin{teo}
If a function $f$ is non-sparse, then $d(M_f)\ge n+1$. If, in 
addition, $f$ satisfies an equation $\varphi (f)=0$, $0\ne \varphi \in 
\AG_{n+1}$, then $f$ is quasi-holonomic.
\end{teo}

\medskip
As in the classical situation, one can construct quasi-holonomic 
functions by addition.

\medskip
\begin{prop}
If the functions $f,g\in \FF_{n+1}$ are quasi-holonomic, and $f+g$ is 
non-sparse, then $f+g$ is quasi-holonomic.
\end{prop}

\medskip
{\it Proof}. Consider the $\AG_{n+1}$-module 
$M_2=(\AG_{n+1}f)\oplus (\AG_{n+1}g)$. Since $f$ and $g$ are 
both quasi-holonomic, we have $d(M_2)=n+1$. Next, let $N_2$ be a 
submodule of $M_2$ consisting of such pairs $(\varphi (f),\varphi 
(g))$ that $\varphi (f)+\varphi (g)=0$. Then $d(M_2)=\max \{ 
d(N_2),d(M_2/N_2)\}$, so that $d(M_2/N_2)\le n+1$.

On the other hand, we have an injective mapping $\AG_{n+1}(f+g)\to 
M_2/N_2$, which maps $\varphi (f+g)$ to the image of $(\varphi (f),\varphi 
(g))$ in $M_2/N_2$. Therefore $d(\AG_{n+1}(f+g))\le d(M_2/N_2)\le 
n+1$. It remains to use Theorem 3. $\qquad \blacksquare$

\medskip
{\bf 4.3.} We use Theorem 3 to prove that the functions (5) 
obtained via the sequence-to-function transform ($*$) or its 
multi-index generalizations, from some well-known sequences of 
polynomials over $K$ are quasi-holonomic.

a) {\it The Carlitz polynomials}. The sequence
$$
f_k(s)=D_k^{-1}\prod \limits _{\genfrac{}{}{0pt}{1}{m\in 
\F [x]}{\deg m<k}}(s-m) \quad (k\ge 1),\quad f_0(s)=s,
$$
of normalized Carlitz polynomials forms an orthonormal basis of 
the space of all $\F$-linear continuous functions on the ring of 
integers of the field $K$. Its transform ($*$), the function
\begin{equation}
C_s(t)=\sum\limits_{k=0}^\infty f_k(s)t^{q^k}
\end{equation}
called the {\it Carlitz module}, is one of the main objects of the 
function field arithmetic \cite{G2,Th3}.

It is known \cite{Carl,G1} that
$$
f_k(s)=\sum\limits_{i=0}^k\frac{(-1)^{k-i}}{D_iL_{k-i}^{q^i}}s^{q^i}
$$
where $L_i=[i][i-1]\ldots [1]$ ($i\ge 1$), $L_0=1$. By (1), we 
have
\begin{equation}
|D_i|=q^{-\frac{q^i-1}{q-1}},\quad |L_i|=q^{-i},
\end{equation}
so that
$$
\left| D_iL_{k-i}^{q^i}\right| =q^{-\left( 
\frac{q^i-1}{q-1}+(k-i)q^i\right) },\quad 0\le i\le k.
$$

For large values of $k$, an elementary investigation of the 
function $z\mapsto (k-z)q^z$, $z\le k$, shows that
$$
\max\limits_{0\le i\le k}(k-i)q^i\le \alpha q^k,\quad \alpha >0,
$$
so that
$$
|f_k(s)|\le q^{\alpha q^k}
$$
for all $s\in \K$ with $|s|\le q^{-1}$. Therefore the series (21) 
converges for small $|t|$, so that the Carlitz module function 
belongs to $\FF_2$.

Since $d_sf_i=f_{i-1}$ for $i\ge 1$, and $d_sf_0=0$ \cite{G1}, we 
see that $d_sC_s(t)=C_s(t)$. Clearly, the function $C_s(t)$ is 
non-sparse. Therefore the Carlitz module function is quasi-holonomic, 
jointly in both its variables.

\medskip
b) {\it Thakur's hypergeometric polynomials}. We consider the 
polynomial case of Thakur's hypergeometric function 
\cite{Th1,Th2,Th3}:
\begin{equation}
{}_lF_\lambda (-a_1,\ldots ,-a_l;-b_1,\ldots ,-b_\lambda 
;z)=\sum\limits_m\frac{(-a_1)_m\ldots (-a_l)_m}{(-b_1)_m\ldots (-b_\lambda 
)_mD_m}z^{q^m}
\end{equation}
where $a_1,\ldots ,a_l,b_1,\ldots ,b_\lambda \in\mathbb Z_+$,
\begin{equation}
(-a)_m=\begin{cases}
(-1)^{a-m}L_{a-m}^{-q^m}, & \text{if $m\le a$},\\
0, & \text{if $m>a$},\end{cases},\quad a\in \mathbb Z_+.
\end{equation}

It is seen from (24) that the terms in (23), which make sense and 
do not vanish, are those with $m\le \min (a_1,\ldots ,a_l,b_1,\ldots 
,b_\lambda )$. Let
\begin{multline}
f(s,t_1,\ldots ,t_l,u_1,\ldots ,u_\lambda )\\ =
\sum\limits_{k_1=0}^\infty \ldots \sum\limits_{k_l=0}^\infty
\sum\limits_{\nu_1=0}^\infty \ldots \sum\limits_{\nu_\lambda =0}^\infty
{}_lF_\lambda (-k_1,\ldots ,-k_l;-\nu_1,\ldots ,-\nu_\lambda ;s)
t_1^{q^{k_1}}\ldots t_l^{q^{k_l}}u_1^{q^{\nu_1}}\ldots 
u_\lambda^{q^{\nu_\lambda}}.
\end{multline}
We prove as above that all the series in (25) converge near the 
origin. Thus, $f\in \FF_{l+\lambda +1}$.

It is known (\cite{Th3}, Sect. 6.5) that
\begin{equation}
d_s{}_lF_\lambda (-k_1,\ldots ,-k_l;-\nu_1,\ldots ,-\nu_\lambda 
;s)={}_lF_\lambda (-k_1+1,\ldots ,-k_l+1;-\nu_1+1,\ldots ,-\nu_\lambda +1;s)
\end{equation}
if all the parameters $k_1,\ldots ,k_l,\nu_1,\ldots ,\nu_\lambda$ 
are different from zero. If at least one of them is equal to zero, 
then the left-hand side of (26) equals zero. This property implies 
the identity $d_sf=f$, the same as that for the Carlitz module 
function. Since $f$ is non-sparse, it is quasi-holonomic.

In the next section we will see that the $K$-binomial coefficients 
(3) correspond to a quasi-holonomic function satisfying a more 
complicated equation containing also the operator $\Delta_t$.

\section{$K$-Binomial Coefficients}

{\bf 5.1.} Let us consider the $K$-binomial coefficients (3). It 
follows from (22) that 
$$
\left| \binom{k}{m}_K\right| =1,\quad 0\le m\le k.
$$
Since $\binom{k}{m}_K\in \F (x)$, it is natural to consider also 
other places of $\F (x)$, that is other non-equivalent absolute 
values on $\F (x)$. It is well known (\cite{Weil}, Sect. 3.1) that 
they are parametrized by monic irreducible polynomials $\pi \in \F 
[x]$. The absolute value $|t|_\pi$, $t\in \F (x)$, is defined as 
follows. We write $t=\pi^\nu \alpha /\alpha'$ where $m\in \mathbb 
Z$, $\alpha ,\alpha'\in \F [x]$, and $\pi$ does not divide $\alpha 
,\alpha'$. Then $|t|_\pi =|\pi |_\pi^\nu$, $|\pi |_\pi 
=q^{-\delta}$ where $\delta =\deg \pi$; as usual, $|0|_\pi =0$. 
The absolute value $|\cdot |$ used elsewhere in this paper 
corresponds to $\pi (x)=x$.

\medskip
\begin{prop}
For any  monic irreducible polynomial $\pi \in \F[x]$, 
the $K$-binomial coefficients (3) satisfy the inequality
$$
\left| \binom{k}{m}_K\right|_\pi \le 1,\quad 0\le m\le k.
$$
\end{prop}

\medskip
{\it Proof}. First we compute $|D_m|_\pi$. It follows from Lemma 
2.13 of \cite{LN} that
$$
|[i]|_\pi =\begin{cases}
q^{-\delta }, & \text{if $\delta$ divides $i$},\\
1, & \text{otherwise}.\end{cases}
$$
Writing $m=j\delta +i$, with $i,j\in \mathbb Z_+$, $0\le 
i<\delta$, we find that
\begin{multline*}
|D_m|_\pi =|[j\delta ]|_\pi^{q^i}|[(j-1)\delta ]|_\pi^{q^{\delta 
+i}}\ldots |[\delta ]|_\pi^{q^{(j-1)\delta +i}}=\left\{ q^{-\delta 
}\cdot \left( q^{-\delta}\right)^{q^\delta }\cdot \ldots \cdot 
\left( q^{-\delta}\right)^{q^{(j-1)\delta }}\right\}^{q^i} \\
=\left\{ \left( q^{-\delta}\right)^{1+q^\delta +\cdots 
+q^{(j-1)\delta}}\right\}^{q^i}=q^{-\delta q^i\frac{q^{j\delta }
-1}{q^\delta -1}}.
\end{multline*}

Similarly we can write $k-m=\varkappa \delta +\lambda$, with 
$\varkappa ,\lambda \in \mathbb Z_+$, $0\le \lambda <\delta$, and 
get that
$$
|D_{k-m}|=q^{-\delta q^\lambda \frac{q^{\varkappa \delta }-1}{q^\delta -1}}.
$$

If $i+\lambda <\delta$, then we obtain a similar representation 
for $k$ simply by adding those for $m$ and $k-m$, so that
\begin{multline*}
\log_q\left| \binom{k}{m}_K\right|_\pi =-\frac{\delta}{q^\delta 
-1}\left\{ q^{i+\lambda }\left( q^{(j+\varkappa )\delta }-1\right) 
-q^i\left( q^{j\delta }-1\right) -q^\lambda \left( q^{\varkappa \delta 
}-1\right) q^{j\delta +i}\right\} \\
=-\frac{\delta}{q^\delta -1}q^i\left( 1+q^{\lambda +j\delta 
}-q^\lambda -q^{j\delta }\right) =-\frac{\delta}{q^\delta 
-1}q^i\left( q^\lambda -1\right) \left( q^{j\delta }-1\right) \le 0.
\end{multline*}

If $i+\lambda \ge \delta$, then $k=(j+\varkappa +1)\delta +\nu$ 
where $0\le \nu =i+\lambda -\delta <\delta$. In this case
\begin{multline*}
\log_q\left| \binom{k}{m}_K\right|_\pi =-\frac{\delta}{q^\delta 
-1}\left\{ q^\nu \left( q^{(j+\varkappa +1)\delta }-1\right) 
-q^i\left( q^{j\delta }-1\right) -q^\lambda \left( q^{\varkappa \delta 
}-1\right) q^{j\delta +i}\right\} \\
=-\frac{\delta}{q^\delta -1}\left( q^i+q^{\lambda +j\delta +i
}-q^{i+j\delta }-q^\nu \right) <0,
\end{multline*}
since $\nu <i+\lambda$. $\qquad \blacksquare$

\medskip
Below we will use only the valuation with $\pi (x)=x$, that is, as 
above, consider the field $K$.

\medskip
{\bf 5.2.} Let us derive, for the $K$-binomial coefficients 
(3), analogs of the classical Pascal and Vandermonde identities.

\medskip
\begin{prop}
The identity
\begin{equation}
\binom{k}{m}_K=\binom{k-1}{m-1}_K^q+\binom{k-1}{m}_K^qD_m^{q-1}
\end{equation}
holds, if $0\le m\le k$ and it is assumed that 
$\dbinom{k}{-1}_K=\dbinom{k-1}{k}_K=0$.
\end{prop}

\medskip
{\it Proof}. Let $e_m(t)=D_mf_m(t)$ be the ``non-normalized'' 
Carlitz polynomials. They satisfy the main $K$-binomial identity 
\cite{Carl,K5}
\begin{equation}
e_k(st)=\sum\limits_{m=0}^k\binom{k}{m}_Ke_m(s)\left\{ 
e_{k-m}(t)\right\}^{q^m},
\end{equation}
which holds, for example, for any $s,t\in \F [x]$.

It is known \cite{Carl,G1} that
\begin{equation}
e_k=e_{k-1}^q-D_{k-1}^{q-1}e_{k-1}.
\end{equation}
Let us rewrite the left-hand side of (28) in accordance with (29), 
and apply to each term the identity (28) with $k-1$ substituted 
for $k$. We have
$$
e_{k-1}^q(st)=\sum\limits_{i=0}^{k-1}\binom{k-1}{i}_K^qe_i^q(s) 
e_{k-i-1}^{q^{i+1}}(t).
$$
By (29), $e_i^q=e_{i+1}+D_i^{q-1}e_i$, 
$e_{k-i-1}^q=e_{k-i}+D_{k-i-1}^{q-1}e_{k-i-1}$, whence
\begin{multline*}
e_{k-1}^q(st)
=\sum\limits_{j=1}^k\binom{k-1}{j-1}_K^qe_j(s)e_{k-j}^{q^j}(t)
+\sum\limits_{i=0}^{k-1}\binom{k-1}{i}_K^q
D_i^{q-1}e_i(s)e_{k-i}^{q^i}(t)\\
+\sum\limits_{i=0}^{k-1}\binom{k-1}{i}_K^q
D_i^{q-1}D_{k-i-1}^{q^i(q-1)}e_i(s)e_{k-i-1}^{q^i}(t).
\end{multline*}

Note that
\begin{equation}
\binom{k-1}{i}_K^qD_i^{q-1}D_{k-i-1}^{q^i(q-1)}=D_{k-1}^{q-1}\binom{k-1}{i}_K.
\end{equation}
Indeed, the left-hand side of (30) equals
$$
\frac{D_{k-1}^q}{D_i^qD_{k-i-1}^{q^{i+1}}}D_i^{q-1}D_{k-i-1}^{q^{i+1}-q^i}=
\frac{D_{k-1}}{D_iD_{k-i-1}^{q^i}}D_{k-1}^{q-1}
$$
and coincides with the right-hand side. Therefore the last sum in 
the expression for $e_{k-1}^q(st)$ equals
$$
D_{k-1}^{q-1}\sum\limits_{i=0}^{k-1}\binom{k-1}{i}_Ke_i(s)e_{k-i-1}^{q^i}(t)
=D_{k-1}^{q-1}e_{k-1}(st).
$$
Using (29) again we find that
$$
e_k(st)=\sum\limits_{i=0}^k\binom{k-1}{i-1}_K^qe_i(s)e_{k-i}^{q^i}(t)
+\sum\limits_{i=0}^k\binom{k-1}{i}_K^qD_i^{q-1}e_i(s)e_{k-i}^{q^i}(t),
$$
and the comparison with (28) yields
$$
\sum\limits_{m=0}^k\left\{ \binom{k}{m}_K-\binom{k-1}{m-1}_K^q-
\binom{k-1}{m}_K^qD_m^{q-1}\right\} e_m(s)e_{k-m}^{q^m}(t)=0
$$
for any $s,t$. 

Since the Carlitz polynomials are linearly independent, we obtain 
that
$$
\left\{ \binom{k}{m}_K-\binom{k-1}{m-1}_K^q-
\binom{k-1}{m}_K^qD_m^{q-1}\right\} e_{k-m}^{q^m}(t)=0
$$
for any $t$, and it remains to note that $e_{k-m}(t)\ne 0$ if 
$t\in \F [x]$, $\deg t\ge k$, by the definition of the Carlitz 
polynomials. $\qquad \blacksquare$

\medskip
More generally, we have the following Vandermonde-type identity. 
Let $k,m$ be integers, $0\le m\le k$.

\begin{prop}
Define $c_{li}^{(m)}\in K$ by the recurrent relation
\begin{equation}
c_{l+1,i}^{(m)}=c_{l,i-1}^{(m)}+c_{li}^{(m)}D_{m-i}^{q-1}
\end{equation}
and the initial conditions $c_{li}^{(m)}=0$ for $i<0$ and $i>l$, 
$c_{00}^{(m)}=1$. Then, for any $l\le m$,
\begin{equation} 
\binom{k}{m}_K=\sum\limits_{i=0}^lc_{li}^{(m)}\binom{k-l}{m-i}_K^{q^l}.
\end{equation}
\end{prop}

\medskip
{\it Proof}. The identity (32) is trivial for $l=0$. Suppose it 
has been proved for some $l$. Let us transform the right-hand side 
of (32) using the identity (27). Then we have
\begin{multline*}
\binom{k}{m}_K=\sum\limits_{i=0}^lc_{li}^{(m)}\binom{k-l-1}{m-i-1}_K^{q^{l+1}}
+\sum\limits_{i=0}^lc_{li}^{(m)}\binom{k-l-1}{m-i}_K^{q^{l+1}}D_{m-i}^{q-1}\\
=\sum\limits_{j=1}^{l+1}c_{l,j-1}^{(m)}\binom{k-l-1}{m-j}_K^{q^{l+1}}
+\sum\limits_{i=0}^lc_{li}^{(m)}\binom{k-l-1}{m-i}_K^{q^{l+1}}D_{m-i}^{q-1}.
\end{multline*}
Since we assume that $c_{l,-1}^{(m)}=c_{l,l+1}^{(m)}=0$, the 
summation in both the above sums can be performed from 0 to $l+1$. 
Using (31) we obtain the required identity (32) with $l+1$ 
substituted for $l$. $\qquad \blacksquare$

\medskip
{\bf 5.3.} Now we consider a function $f\in \FF_2$ associated with 
the $K$-binomial coefficients, that is
\begin{equation} 
f(s,t)=\sum\limits_{k=0}^\infty 
\sum\limits_{m=0}^k\binom{k}{m}_Ks^{q^m}t^{q^k}.
\end{equation}
Obviously, $f$ is non-sparse.

\medskip
\begin{prop}
The function (33) satisfies the equation
\begin{equation} 
d_sf(s,t)=\Delta_tf(s,t)+[1]^{1/q}f(s,t),
\end{equation}
so that $f$ is quasi-holonomic.
\end{prop}

\medskip
{\it Proof}. Let us compute $d_sf$. We have
$$
d_sf(s,t)=\sum\limits_{k=1}^\infty 
\sum\limits_{m=1}^k\binom{k}{m}_K^{1/q}[m]^{1/q}s^{q^{m-1}}t^{q^{k-1}}
=\sum\limits_{\nu =0}^\infty 
\sum\limits_{\mu =0}^\nu \binom{\nu +1}{\mu +1}_K^{1/q}[\mu 
+1]^{1/q}s^{q^\mu}t^{q^\nu}.
$$
Using Proposition 5 we find that $d_sf=\Sigma_1+\Sigma_2$ where
$$
\Sigma_1=\sum\limits_{\nu =0}^\infty 
\sum\limits_{\mu =0}^\nu \binom{\nu}{\mu}_K[\mu 
+1]^{1/q}s^{q^\mu}t^{q^\nu},
$$
$$
\Sigma_2=\sum\limits_{\nu =0}^\infty 
\sum\limits_{\mu =0}^\nu \binom{\nu}{\mu +1}_K[\mu 
+1]^{1/q}D_{\mu +1}^{1-q^{-1}}s^{q^\mu}t^{q^\nu}.
$$

Note that
$$
[\mu +1]^{1/q}=\left( x^{q^{\mu +1}}-x\right)^{1/q}=\left( 
x^{q^\mu }-x\right) +\left( x^q-x\right)^{1/q}=[\mu ]+[1]^{1/q},
$$
so that
\begin{equation}
\Sigma_1=\sum\limits_{\nu =0}^\infty 
\sum\limits_{\mu =0}^\nu \binom{\nu}{\mu}_K[\mu ]
s^{q^\mu}t^{q^\nu}+[1]^{1/q}f(s,t).
\end{equation}

Next, we have
$$
\binom{\nu}{\mu +1}_K[\mu +1]^{1/q}D_{\mu 
+1}^{1-q^{-1}}=\frac{D_\nu}{D_{\mu +1}D_{\nu -\mu -1}^{q^{\mu 
+1}}}D_{\mu +1}\left( \frac{[\mu +1]}{D_{\mu +1}}\right)^{1/q}
=\frac{D_\nu}{D_\mu D_{\nu -\mu -1}^{q^{\mu +1}}},
$$
and also
$$
D_{\nu -\mu -1}^q=\frac{1}{[\nu -\mu ]}[\nu -\mu ]D_{\nu -\mu -1}^q
=\frac{D_{\nu -\mu}}{[\nu -\mu ]},
$$
whence
$$
D_{\nu -\mu -1}^{q^{\mu +1}}=\frac{D_{\nu -\mu }^{q^\mu}}{[\nu -\mu 
]^{q^\mu }}.
$$
Therefore
$$
\Sigma_2=\sum\limits_{\nu =0}^\infty 
\sum\limits_{\mu =0}^\nu \binom{\nu}{\mu}_K[\nu -\mu 
]^{q^\mu}s^{q^\mu}t^{q^\nu}.
$$
As above, $[\nu -\mu ]^{q^\mu }=\left( x^{q^{\nu 
-\mu}}-x\right)^{q^\mu}=[\nu ]-[\mu ]$, so that
$$
\Sigma_2=\sum\limits_{\nu =0}^\infty 
\sum\limits_{\mu =0}^\nu \binom{\nu}{\mu}_K([\nu ]-[\mu 
])s^{q^\mu}t^{q^\nu}.
$$
Together with (35), this implies (34). $\qquad \blacksquare$

\end{document}